\documentclass[11pt, a4paper]{article}
\usepackage{indentfirst}
\usepackage{epsfig}
\usepackage{subfigure}
\usepackage{amssymb}
\usepackage{amsthm}
\usepackage{mathrsfs}
\usepackage{hhline}
\usepackage{longtable}
\usepackage{amsmath}
\usepackage{array}
\usepackage{float}
\usepackage{multirow}
\usepackage{cite}
\usepackage{enumerate}
\usepackage{xcolor}
\usepackage{caption}
\numberwithin{equation}{section}

\newtheorem{theorem}{Theorem}[section]
\newtheorem{corollary}[theorem]{Corollary}
\newtheorem{lemma}[theorem]{Lemma}

\newtheorem{problem}[theorem]{Problem}

\newtheorem{example}{Example}[section]
\newtheorem{conjecture}{Conjecture}[section]

\newtheorem{defi}{Definition}[section]

\usepackage{xcolor}

\definecolor{vividviolet}{rgb}{0.62, 0.0, 1.0}

\def\Z{\mathbb{Z}}

\newcolumntype{"}{@{\hskip\tabcolsep\vrule width 1pt\hskip\tabcolsep}}
\def\ms{\medskip}
\def\nt{\noindent}
\def\D{\Delta}
\def\d{\delta}

\graphicspath{%
    {converted_graphics/}
    {/}
}
\begin{document}
\baselineskip18truept
\normalsize
\begin{center}
{\mathversion{bold}\Large \bf Complete characterization of graphs with local total antimagic chromatic number 3 }

\bigskip
{\large Gee-Choon Lau}

\medskip

\emph{77D, Jalan Subuh,}\\
\emph{85000, Johor, Malaysia.}\\
\emph{geeclau@yahoo.com}\\

\ms\ms

\emph{Dedicated to Prof. Wai Chee Shiu on the occasion of his 66th birthday}
\ms

\end{center}
%

\begin{abstract}
A total labeling of a graph $G = (V, E)$ is said to be local total antimagic if it is a bijection $f: V\cup E \to\{1,\ldots ,|V|+|E|\}$ such that adjacent vertices, adjacent edges, and incident vertex and edge have distinct induced weights where the induced weight of a vertex $v$ is $w_f(v) = \sum f(e)$ with $e$ ranging over all the edges incident to $v$, and the induced weight of an edge $uv$ is $w_f(uv) = f(u) + f(v)$.  The local total antimagic chromatic number of $G$, denoted by $\chi_{lt}(G)$, is the minimum number of distinct induced vertex and edge weights over all local total antimagic labelings of $G$.  In this paper, we first obtained general lower and upper bounds for $\chi_{lt}(G)$ and sufficient conditions to construct a graph $H$ with $k$ pendant edges and $\chi_{lt}(H) \in\{\D(H)+1, k+1\}$. We then completely characterized graphs $G$ with $\chi_{lt}(G)=3$. Many families of (disconnected) graphs $H$ with $k$ pendant edges and $\chi_{lt}(H) \in\{\D(H)+1, k+1\}$ are also obtained.  

\ms

\end{abstract}

\baselineskip18truept
\normalsize


\section{Introduction}

\nt Let $G=(V, E)$ be a simple and loopless graph of order $p$ and size $q$. For integers $a < b$, let $[a,b] = \{n\in\Z\;|\; a\le n\le b\}$. If no adjacent vertices of $G$ are assigned the same color, then $G$ has a proper (vertex) coloring. The smallest required number of color(s) is the chromatic number of $G$, denoted $\chi(G)$. The most famous problem on chromatic number is the 4-Color Conjecture which states that every planar graph $G$ has $\chi(G) = 4$. Interested readers may refer to~\cite{AppelHaken} for the proof by Appel and Haken. If no two adjacent edges of $G$ are assigned the same color, then $G$ has a proper edge coloring. The smallest required number of color(s) is the chromatic index of $G$, denoted $\chi'(G)$.  Interested readers may refer to~\cite{ToftWilson} for an interesting brief history. A total coloring of $G$ is a coloring of the vertices and edges of $G$ such that $x$ and $y$ are assigned distinct colors whenever $x$ and $y$ are adjacent vertices or edges, or incident vertex and edge. The smallest number of color(s) required for a total coloring is the total chromatic number of $G$, denoted $\chi_t(G)$. The conjecture that $\D(G)+1 \le \chi_t(G) \le \D(G)+2$ remained unsolved (see~\cite{Behzad, G+N+S}) where $\Delta(G)$ is the maximum degree of $G$. 

\ms\nt A bijection $f : E  \to [1,q]$ is a local antimagic labeling if $f^+(u) \ne f^+(v)$ for every edge $uv$ of $G$ where $f^+(u)$ is the sum of all the incident edge label(s) of $u$ under $f$. Let $f^+(u)$ be the color of $u$. Clearly, $f^+$ is a proper vertex coloring of $G$ induced by $f$. The smallest number of vertex colors taken over all local antimagic labeling $f$ is called the local antimagic chromatic number of $G$, denoted $\chi_{la}(G)$ (see~\cite{Arumugam, Haslegrave,LNS, LSN, LS-AMH1, LS-AMH2} for the many results available). 

\ms\nt A bijection $g: V \to [1,p]$ is a local edge antimagic labeling if $g_{le}(e_1)\ne g_{le}(e_2)$ for every two adjacent edges $e_1$ and $e_2$ of $G$ where $g_{le}(uv) = g(u)+g(v)$. Let $g_{le}(e)$ be the induced edge color of $e$ under $g$. Clearly, $g_{le}$ is a proper edge coloring of $G$ induced by $g$.  The smallest number of edge color(s) over all local edge antimagic labeling $g$ is called the local edge antimagic chromatic number of $G$, denoted $\chi_{lea}(G)$~\cite{Agustin}.

\ms\nt A bijection $f:V\cup E\to [1,p+q]$ is a local antimagic total labeling if $w^+(u)\ne w^+(v)$ for every two adjacent vertices $u$ and $v$ where $w^+(u) = f(u) + \sum f(e)$ over all edge(s) $e$ incident to $u$. Let $w^+(u)$ be the induced vertex color of $u$ under $f$. Clearly, $w^+$ is a proper vertex coloring of $G$ induced by $f$. The smallest number of vertex colors over all local antimagic total labeling $f$ is called the local antimagic total chromatic number of $G$, denoted $\chi_{lat}(G)$. In~\cite{LS-OM,LSS-OM}, the authors proved that every graph is local antimagic total and obtained the exact $\chi_{lat}(G)$ for many families of graphs $G$.

\ms\nt Motivated by the concepts of local (edge) antimagic labeling and total labeling of $G$, Sandhiya and Nalliah~\cite{S+N} defined the concept of local total antimagic labeling.  A bijection $f: V\cup E\rightarrow [1, p+q]$ such that the weight of a vertex $u$ is $w_f(u) = \sum f(e)$ over all edge(s) $e$ incident to $u$, and the weight of an edge $e=uv$ is  $w_f(e) = f(u) + f(v)$, is called a {\it local total antimagic labeling} if (i) every two adjacent vertices and adjacent edges have distinct weights, and (ii) every vertex and its incident edge(s) have distinct weights.  The mapping $w_f$ (or $w$ if no ambiguity) is called a \textit{local total antimagic labeling of $G$ induced by $f$}, and the weights assigned to vertices and edges are called \textit{induced total colors} under $f$. 

\ms\nt The \textit{local total color number} of a local total antimagic labeling $f$ is the number of distinct induced total colors under $f$, denoted by $w(f)$.  Moreover, $f$ is called a {\it local total antimagic $w(f)$-coloring} and $G$ is {\it local total antimagic $w(f)$-colorable}. The {\it local total antimagic chromatic number} $\chi_{lt}(G)$ is defined to be the minimum number of colors taken over all local total antimagic colorings of $G$ induced by local total antimagic labelings of $G$. Let $G+H$ and $mG$ denote the disjoint union of graphs $G$ and $H$, and $m$ copies of $G$, respectively. By definition, any graph with an isolated vertex or a $K_2$ component does not admit a local total antimagic labeling. In this paper, we only consider graph $G$ of order $p\ge 3$ and size $q\ge 2$. Thus, $G$ contains neither isolated vertex nor $K_2$ component.   
 
\ms\nt In Section~\ref{sec-SC}, we obtained various sufficient conditions on bounds of $\chi_{lt}(G)$. In Section~\ref{sec-chilt3}, we give complete characterization of graphs $G$ with $\chi_{lt}(G) = 3$. In Section~\ref{sec-pend+1}, we first obtained various families of disconnected graph with $k$ pendant edges and $\chi_{lt}(G) = k+1$. Using the results in Section~\ref{sec-chilt3}, we then obtained various new families of graphs $H$, constructed from $G$, having $s$ pendant edges and $\chi_{lt}(G) = s+1$. Open problems and conjectures are given in Section~\ref{sec-con}.

\section{Sufficient conditions}\label{sec-SC}

\ms\nt In~\cite{S+N}, the authors proved that a graph $G$ with $k\ge 1$ pendant edge(s) has $\chi_{lt}(G)\ge k+1$. We now give a necessary condition for equality to hold.
 
\begin{lemma}\label{lem-pendant} Suppose $G$ has $k\ge 1$ pendant edges, then $\chi_{lt}(G)\ge k+1$. If $\chi_{lt}(G)=k+1$, then $p+q$ is assigned to a pendant vertex or a pendant edge. \end{lemma}

\begin{proof} Assume that $\chi_{lt}(G) = k+1$. Suppose $p+q$ is assigned to a non-pendant vertex, say $x$. Now, all the $k$ pendant vertices have distinct weights at most $p+q-1$. Moreover, $x$ is incident to two edges, say $e$ and $e'$. Clearly, $w(e), w(e')$ are distinct weights at least $p+q+1$. Suppose $p+q$ is assigned to a non-pendant edge, say $xy$. Now, all the $k$ pendant vertices have distinct weights at most $p+q-1$ and $w(x),w(y)$ are distinct weights at least $p+q+1$. In both cases, there are at least $k+2$ distinct weight so that $\chi_{lt}(G)\ge k+2$. This contradicts $\chi_{lt}(G) = k+1$. \end{proof}

\begin{corollary}\label{cor-pendant} Suppose $G$ is a graph with $k\ge 1$ pendant edges, then $\chi_{lt}(G)\ge \max\{\Delta(G)+1, k+1\}\ge 3$. \end{corollary} 

\begin{proof} By definition and Lemma~\ref{lem-pendant}, $\chi_{lt}(G)\ge \chi_t(G)\ge \max\{\Delta(G)+1,k+1\}$. Thus, if $\chi_{lt}(G)=2$, then $\Delta(G)=1$ or $k=1$. If the former holds, $G$ is a 1-regular graph that does not admit a local total antimagic labeling. If the later holds, $\D(G)\ge 3$ so that $\chi_{lt}(G) \ge \chi_t(G)\ge 4$, a constradiction. Consequently, $\chi_{lt}(G)\ge 3$ if $G$ admits a local total antimagic labeling. \end{proof}

\nt In~\cite{S+N}, the authors also proved that if $G$ is an $r$-regular graph, $r\ge 2$, that admits a $\chi_{la}(G)$-labeling and a $\chi_{lea}(G)$-labeling, then $\chi_{lt}(G) \le \chi_{la}(G) + \chi_{lea}(G)$. We now give a more general result.

\begin{theorem}\label{thm-chilt} Suppose $G$ admits a $\chi_{lea}(G)$-labeling. If G has (i) $\d(G)\ge 2$; or else (ii) $\d(G) = 1$ so that $G$ has $k\ge 1$ pendant vertices, $v_p$ non-pendant vertices and $e_p$ non-pendant edges with $e_p > v_p + k - 2$ that admits a $\chi_{la}(G)$-labeling that assigns the non-pendant edges by integers in $[1,e_p]$, then $G$ admits a local total antimagic labeling. Moreover, $\chi_{lt}(G) \le \chi_{la}(G) + \chi_{lea}(G)$. 
\end{theorem} 

\begin{proof} Suppose $G$ has order $p$ and size $q$. By definition, $G$ has no $K_2$ components.  Let $f : E \to [1,q]$ (and $h : V \to [1,p]$) be a $\chi_{la}(G)$- (and $\chi_{lea}(G)$-) labeling of $G$. Define a bijective total labeling $g : V \cup E \to [1,p+q]$ such that $g(u) = h(u)$ for each vertex $u$ of $G$, and $g(e) = f(e) + p$ for each edge $e$ of $G$. Therefore, $\{g(u)\} = [1,p]$ and $\{g(e)\} = [p+1,p+q]$. Clearly, $w_g(u)= w_g(v)$ if and only if $f^+(u)= f^+(v)$, and $w_g(e_1)= w_g(e_2)$ if and only if  $h_{le}(e_1)= h_{le}(e_2)$ for all $u,v\in V(G)$ and $e_1,e_2\in E(G)$.

\ms\nt {\bf (i).} $\d\ge 2$. In this case, $G$ has no pendant edges. Thus, every edge weight under $g$ is at most $2p-1$ and every vertex weight under $g$ is at least $2p+3$. Therefore, every edge weight is less than every vertex weight. Consequently, $g$ is a local total antimagic labeling that induces $\chi_{la}(G)+\chi_{lea}(G)$ distinct weights. 

\ms\nt {\bf (ii).} $\d=1$. Thus, $G$ has order $p = v_p+k$ and size $q = e_p+k$. By the given assumption, the $\chi_{la}(G)$-labeling $f$ of $G$ assigns each non-pendant edge with an integer in $[1, e_p]$. In this case, under $g$, every vertex is assigned an integer in $[1,v_ p+k]$, every non-pendant edge is assigned an integer in $[v_p+k+1, v_p+k+e_p]$ while every pendant edge is assigned an integer in  $[v_p+k+e_p+1, v_p+2k+e_p]$. Now, every edge weight is at most $2(v_p+k)-1$. Every pendant vertex weight is at least $v_p+k+e_p+1 > 2(v_p+k) - 1$ and every non pendant vertex weight is at least $2(v_p+k)+3$.  Therefore, every edge weight is less than every vertex weight. Consequently, $g$ is a local total antimagic labeling that induces $\chi_{la}(G) + \chi_{lea}(G)$ distinct weights. 
\end{proof}

\nt Suppose $G$ is a graph with $\chi_{lt}(G) = \Delta(G)+1$. We now give sufficient condition to contruct new graph $H$ with $k\ge 2$ pendants edges from $G$  to have $\chi_{lt}(H)=\Delta(H)+1$.

\begin{theorem}\label{thm-HfromG} For $k,s\ge1, ks\ge 2$, let $G$ be a graph of order $p$ and size $q$ with $d\ge 0$ pendant edges such that
\begin{enumerate}[(a)]
\item $\chi_{lt}(G) = t = \Delta(G)+1$ and  the corresponding local total antimagic labeling $f$ assigns $k$ to a maximum degree vertex $v$ of $G$; 
\item for $k\ge 1, s\ge 1$, $w_f(v)+\sum^s_{j=1}\sum^k_{i=1}  (p+q+2jk+1-i)\ne w_f(x)$, for each vertex $x$ adjacent to $v$; 
\item $v$ is the only element with weight $w_f(v)$, or else there is a non-pendant vertex element in $V(G)\cup E(G)$ with weight in $\{p+q+(2j-1)k+i\,|\,1\le j\le s, 1\le i\le k\}$.
\end{enumerate}
Suppose $G_v(k,s)$ is obtained from $G$ by attaching $ks\ge 2$ pendant edges to $v$. If Conditions (a), (b) are satisfied, then $\D(G_v(k,s))+1\le \chi_{lt}(G_v(k,s))\le \Delta(G_v(k,s))+2$. The lower bound holds if Condition (c) is also satisfied. \end{theorem}

\begin{proof} Let $G$ be a graph satisfying Condition (a). Clearly, $v$ is not a pendant vertex.  Lemma~\ref{lem-pendant} implies that $\D(G)\ge d$ if $d>0$.  Note that $G_v(k,s)$ has $ks+d$ pendant edges. So, $\D(G_v(k,s)) = \D(G)+ks\ge ks + d$.  By definition, $\chi_{lt}(G_v(k,s))\ge \Delta(G_v(k,s))+1 = \D(G)+ks+1$.

\ms\nt Let the $ks\ge 2$ pendant edges added to $G$ to get $H$ be $e_{j,i} = vx_{j,i}$ for $1\le j\le s, 1\le i\le k$. For simplicity, if $k=1$, let $e_{j,1} = e_j$ and $x_{j,1} = x_j$. Define a total labeling $g : V(G_v(k,s))\cup E(G_v(k,s)) \to [1, p+q + 2ks]$ such that $g(x) = f(x)$ for each $x\in V(G)\cup E(G)$. For the remaining vertices and edges, we do as follows.  
\begin{enumerate}[(i)]
\item For $k=1$, let $g(x_j) = p+q+2j-1$ and $g(e_j) = p+q+2(s-j+1)$. If $s\ge 3$ is odd, then swap the edge labels of $e_{(s+1)/2}$ and $e_{(s-1)/2}$. 
\item For $k\ge 2$, let $g(x_{j,i}) = p+q+ 2(j-1)k+i$, and $g(e_{j,i}) = p+q+2jk+1-i$ for $1\le j\le s, 1\le i\le k$. If $k\ge 3$ is odd, then swap the edge labels of $e_{j,(k+1)/2}$ and $e_{j,(k-1)/2}$ for each $1\le j\le s$. 
\end{enumerate} 
Clearly, for distinct $x, y\in (V(G)\cup E(G)) \setminus\{v\}$, $w_g(x) = w_g(y)$ if and only if $w_f(x) = w_f(y)$. Observe that $\{g(e_{j,i})\mid 1\le j\le s, 1\le i\le k\} = \bigcup^s_{j=1} [p+q+(2j-1)k+1,p+q+2jk] = \{w_g(e_{j,i}) \mid 1\le j\le s, 1\le i\le k\} = \{w_g(x_{j,i}) \mid 1\le j\le s, 1\le i\le k\}$, denoted $W$. 
Clearly, all the elements in $W$ are distinct. Moreover, every edge $e$ of $G$ that is incident to $v$  has weight $w_g(e) = w_f(e) \le 2p+2q-1< w_f(v) +  \sum^s_{j=1}\sum^k_{i=1} (p+q+(2j-1)k+i) = w_g(v)$. So, Condition (b) implies that $g$ is a local total antimagic labeling of $H$ with weights set $W_g = W \cup (W_f\setminus \{w_f(v)\}) \cup \{w_g(v)\}$ that has size at most $ks + t + 1 = \D(G) + ks +2 = \D(G_v(k,s)) + 2$. Thus, $\D(G_v(k,s))+1\le \chi_{lt}(G_v(k,s))\le \Delta(G_v(k,s))+2$.  

\ms\nt Consider Condition (c).  Suppose $v$ is the only element  with weight $w_f(v)$. Now $g$ induces $t-1$ distinct weights among the elements in $(V(G)\cup E(G))\setminus\{v\}$. Therefore, $W_g$ has size $ks +(t-1)+1 = ks + t = \D(G_v(k,s)) + 1$. Otherwise, since there is a non-pendant vertex element in $(V(G)\cup E(G))\setminus\{v\}$ with weight in $\{p+q+(2j-1)k+i\,|\,1\le j\le s, 1\le i\le k\}$, we also have $W_g$ has size $\D(G_v(k,s)) + 1$. This completes the proof. \end{proof}

\nt Suppose $G$ has $d_G\ge 1$ pendant edges and $\chi_{lt}(G) = d_G+1\ge\D(G)+1$. Let $H$ be obtained from $G$ by attaching pendant edges to a vertex of $G$ such that $H$ has $d_H$ pendant edges. We now give sufficient condition for $\chi_{lt}(H) = d_H+1$. Let $\deg_G(v)$ be the degree of vertex $v$ in $G$. 

\begin{theorem}\label{thm-HfromG2} For $k,s\ge 1$, $ks\ge 2$, let $G$ be a graph of order $p$ and size $q$ with $d_G$ pendant edges such that 
\begin{enumerate}[(i)]
\item $\chi_{lt}(G) = d_G+1\ge \D(G)+1\ge 3$ and the corresponding local total antimagic labeling $f$ assigns $k$ to a  vertex $v$ of $G$,
\item if $x$ is a vertex adjacent to $v$, then $w_f(x)\ne w_f(v)+\sum^s_{j=1} \sum^k_{i=1} (p+q+(2j-1)k+i)$, 
\item if $v$ is not the only element with weight $w_f(v)$, then $w$ or else $w_f(v)\in \{g(e_{j,i})\mid 1\le j\le s, 1\le i\le k\}$, where $w$ is the only weight under $f$ which is not a pendant edge label. 
\end{enumerate}
Suppose $G_v(k,s)$ is obtained from $G$ by attaching $ks$ pendant edges to $v$. 
\begin{enumerate}[(a)]
\item If $v$ is a pendant vertex, then $d_H+1=ks+d_G\le \chi_{lt}(G_v(k,s))\le ks+d_G+1=d_H+2$. The lower bound holds if  (i) there is a weight of $G$ under $f$ equal to  $p+q+(2j-1)k+i$ for $1\le j\le s, 1\le i\le k$; or else (ii) $v$ is the only element with weight $w_f(v)$ under $f$. 
\item If $v$ is not a pendant vertex, then  $ \chi_{lt}(G_v(k,s)) = ks+d_G+1$. 
\end{enumerate}
\end{theorem}

\begin{proof} Let the $ks$ vertices added to $G$ to get $H=G_v(k,s)$ be $x_{j,i}$ for $1\le j\le s, 1\le i\le k$ and let $vx_{j,i}$ be $e_{j,i}$. Since $H$ has order $p+ks$ and size $q+ks$, we define a total labeling $g: V(H) \cup E(H) \to [1, p+q+2ks]$ such that $g$ is as defined in the proof of Theorem~\ref{thm-HfromG}. Obviously, every edge incident to $v$ has weight less than $w_g(e_{j,i})$. Moreover, all but one of the weights of $G$ under $f$, say $w$, must be an edge label of $G$ under $f$. Moreover, $w_g(v) = w_f(v) + \sum^s_{j=1}\sum^k_{i=1} (p+q+(2j-1)k+i)$. 

\ms\nt By definition and Lemma~\ref{lem-pendant}, we know $d_G\ge \D(G)\ge 2$.
\begin{enumerate}[(a)]
\item If $v$ is a pendant vertex, then $H$ has $d_H = ks+d_G-1\ge ks+1$ pendant edges with $\deg_H(v) = ks+1$. Since $ks\ge 2$, we also have $ks+d_G-1\ge d_G+1 > \D(G)$ so that $ks+d_G-1\ge \max\{ks+1,\D(G)\}=\D(H)$. By Lemma~\ref{lem-pendant}, $\chi_{lt}(H)\ge ks+d_G = d_H + 1$. Note that $w_g(x) =w_g(y)$ if and only if $w_f(x) = w_f(y)$ for each $x, y \in (V(G)\cup E(G))\setminus\{v\}$. Clearly, $w_f(v)$ is unique. If $w$ does not equal to   $p+q+(2j-1)k+i$ for $1\le j\le s, 1\le i\le k$, then  Condition (ii) implies that $g$ is a local total antimagic labeling that induces $d_G+ks+1$ distinct weights so that $d_G+ks\le \chi_{lt}(H)\le d_G+ks+1$. Otherwise, $g$ induces $d_G+ks$ distinct weights and the lower bound holds.

\item If $v$ is not a pendant vertex, then $\D(G)\ge \deg_G(v) = r\ge 2$ and $H$ has $ks+d_G$ pendant edges with $\deg_H(v) = ks + r$ so that $ks+d_G \ge \max\{ks + r,\D(G)\} = \D(H)$. By Lemma~\ref{lem-pendant}, $\chi_{lt}(H)\ge ks+d_G+1$. Note that $w_g(x) =w_g(y)$ if and only if $w_f(x) = w_f(y)$ for each $x, y \in (V(G)\cup E(G))\setminus\{v\}$. If $v$ is the only element with weight $w_f(v)$, then Condition (ii) implies that $g$ is a local total antimagic labeling that induces $ks+d_G+1$ distinct weights. Otherwise, Condition (iii) implies that $g$ is a local total antimagic chromatic labeling that induces $ks+d_G+1$ distinct weights. Thus,  $\chi_{lt}(H) = ks+d_G+1$.
\end{enumerate}
This completes the proof.
\end{proof}

\begin{theorem}\label{thm-D} Suppose $G$ is a graph of order $p$ and size $q$ with exactly one vertex of maximum degree $\D\ge 3$ which is not adjacent to any pendant vertex and all other vertices of $G$ has degree at most $m < \D$. Moreover, $G$ has $k\ge \D\ge 2$ pendant edges  such that $\D(\D+1) > \max\{m[2(p+q)-m+1], 4(p+q)-2\}$. If $G$ admits a local total antimagic labeling, then $\chi_{lt}(G)\ge k+2$.
\end{theorem}

\begin{proof} Let $f$ be a local total antimagic labeling of $G$. If $p+q$ is assigned to a non-pendant vertex or edge, by Lemma~\ref{lem-pendant}, $\chi_{lt}(G)\ge k+2$. Suppose $p+q$ is assigned to a pendant edge that has a non-pendant end-vertex $x$. Now, all the $k$ induced pendant vertex labels are distinct and at most $p+q$. Suppose  $u$ is the vertex of maximum degree $\D$, then $w(u)\ge \D(\D+1)/2$ and $p+q+1 \le w(x) \le m[2(p+q)-m+1]/2$. By the given hypothesis, $w(u) > w(x) > w(y)$ for every pendant vertex $y$. Thus, $f$ induces at least $k+2$ distinct vertex weights. 

\ms\nt Suppose $p+q$ is assigned to a pendant vertex, then the adjacent pendant edge, say $e$, has $p+q+1\le w(e)\le 2(p+q)-1$. Moreover, all the $k$ induced penant vertex labels are distinct and at most $p+q-1$. By the given hypothesis, $w(u) \ge \D(\D+1)/2 > 2(p+q)-1\ge w(e)$. Thus, $f$ induces at least $k+2$ distinct vertex weights. 
\end{proof}

\nt Let $f_n$ be the fan graph obtained from $n\ge 2$ copies of $K_3$ with a common vertex $c$. Let $f_n(k), k \le 2n-3,$ be obtained from $f_n$ by attaching exactly $k$ pendant edges to every degree 2 vertex of $f_n$.  Now, $f_n(k)$ has order $p= n(2k+2)+1$ and size $q = n(2k+3)$ with exactly one vertex $c$ with maximum degree $\D = 2n$ with is not adjacent to any pendant vertex and all other vertices has degree at most $m = k+2\le 2n-1 < \D$. Thus, $\max\{m[(2(p+q)-m+1], 4(p+q)-2\} = m[2(p+q)-m+1] = 2n(k+2)(4k+5)-(k+2)(k-1)$. Thus, $\D(\D+1) >m[2(p+q)-m+1]$ implies that $2n(2n+1) - 2n(k+2)(4k+5) + (k+2)(k-1) > 0$.   By Lemma~\ref{lem-pendant} and Theorem~\ref{thm-D}, we have the following.   

\begin{corollary}\label{cor-D} For $2n\ge k+3\ge 4$, and $2n(2n+1) - 2n(k+2)(4k+5) + (k+2)(k-1) > 0$, $\chi_{lt}(f_n(k))\ge 2nk+2$. \end{corollary}

\nt Note that the above corollary is always attainable. For example, take $k=1$, we have $n>13$.

\section{Graphs with $\chi_{lt}=3$}\label{sec-chilt3}

\nt By Corollary~\ref{cor-pendant}, $\chi_{lt}(G) = 3$ only if $\D(G) = 2$ possibly with exactly one path component of order at least 3. Let $P_n = u_1u_2\ldots u_n$ be the path of order $n\ge 3$ with $e_i = u_iu_{i+1}$ for $1\le i\le n-1$. The authors in~\cite{S+N} also obtained $\chi_{lt}(P_n)=3$ for  $n=3,6$ and $\chi_{lt}(P_n)=4$ for $n=4$ and odd $n\ge 5$. Moreover, $3\le \chi_{lt}(P_n) \le 5$ for even $n\ge 8$. We first improve the lower bound of the last statement in the following lemma. 

\begin{lemma}\label{lem-Pn} For $n=3,6$, $\chi_{lt}(P_n) = 3$. If $n=4$ or $n\ge 5$ is odd, $\chi_{lt}(P_n) = 4$. Otherwise, $4\le \chi_{lt}(P_n)\le 5$ for even $n\ge 8$. \end{lemma} 

\begin{proof} From the proofs in~\cite[Theorems 2.3 and 2.5]{S+N}, we only need to show that $\chi_{lt}(P_n)\ge 4$ for $n\ge 7$. Suppose $\chi_{lt}(P_n) = 3$ and $f$ is a corresponding local total antimagic 3-coloring. Without loss of generality, we may assume the 3 distinct weights are $a,b,c$ such that $a = w(u_1) = w(e_2) = w(u_4) = w(e_5)$, $b=w(e_1) = w(u_3) = w(e_4) = w(u_6)$ and $c=w(u_2)=w(e_3) = w(u_5) = w(e_6)$. Thus, we have 
\begin{eqnarray*}
a = f(e_1) = f(e_3)+f(e_4)\\
b = f(e_2)+f(e_3) = f(e_5)+f(e_6)\\
c = f(e_1)+f(e_2) = f(e_4)+f(e_5).
\end{eqnarray*} 
 Thus, from $a+b-c$, we get $f(e_3) = f(e_3) + f(e_6)$, a contradiction. Therefore, $\chi_{lt}(P_n)\ge 4$. The lemma holds.
\end{proof}

\nt By Corollary~\ref{cor-pendant}, the proofs in~\cite[Theorems 2.3 and 2.5]{S+N} and the argument of Lemma~\ref{lem-Pn}, we have the following.
\begin{corollary}\label{cor-2regPn} If $\chi_{lt}(G) = 3$, then $G$ is a 2-regular graph or a path $P_n, n = 3,6$, or $H+P_3$ or $H+P_6$, where $H$ is a 2-regular graph.  \end{corollary}


\nt In what follows, if a 2-regular graph has $i$-th component of order $n\ge 3$, then the consecutive vertices and edges are $u_{i,1}, e_{i,1}, u_{i,2}, e_{i,2}, \ldots, u_{i,n}, e_{i,n}$. We first give a family of 2-regular graphs $G$ with $\chi_{lt}(G) = 3$.

\begin{theorem}\label{thm-mC6} For $m\ge 1$, $\chi_{lt}(mC_6) = 3$. \end{theorem}

\begin{proof} Suppose the $i$-th copy of $C_6$ has vertices and edges $u_{i,1},$ $e_{i,1},$ $u_{i,2},$ $e_{i,2},$ $u_{i,3},$ $e_{i,3},$ $u_{i,4},$ $e_{i,4},$ $u_{i,5},$ $e_{i,5},$ $u_{i,6},$ $e_{i,6}$ consecutively for $1\le i\le m$. Since $mC_6$ has $6m$ vertices and $6m$ edges, we define a bijection $f : V(mC_6) \cup E(mC_6) \to [1, 12m]$ such that for $1\le i\le m$,
\begin{enumerate}[(i)]
\item $f(u_{i,1}) = 3i - 2$, $f(u_{i,3}) = 3i-1$, $f(u_{i,5}) = 3i$,
\item $f(u_{i,2}) = 12m + 3 - 3i$, $f(u_{i,4}) = 12m + 1 - 3i$, $f(u_{i,6}) = 12m + 2 - 3i$,
\item $f(e_{i,1}) = 6m - 1 + 3i$, $f(e_{i,3}) = 6m + 3i$, $f(e_{i,5}) = 6m - 2 + 3i$,
\item $f(e_{i,2}) = 6m +1 - 3i$, $f(e_{i,4}) = 6m+2-3i$, $f(e_{i,6}) = 6m+3-3i$. 
\end{enumerate}

\nt Clearly, the weights of $u_{i,1},$ $e_{i,1},$ $u_{i,2},$ $e_{i,2},$ $u_{i,3},$ $e_{i,3},$ $u_{i,4},$ $e_{i,4},$ $u_{i,5},$ $e_{i,5},$ $u_{i,6},$ $e_{i,6}$ are $12m+2$, $12m+1$, $12m$ repeatedly for $1\le i\le m$. Thus, $\chi_{lt}(mC_6)\le 3$. Since $\chi_{lt}(mC_6)\ge \chi_{t}(mC_6) = 3$, the theorem holds.
\end{proof}

\begin{example} The figure below gives the local total antimagic 3-coloring of $2C_6$ as defined above with induced weights $24,25,26$.
\begin{figure}[H]
\begin{center}
\centerline{\epsfig{file=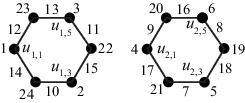, width=5cm}}
\caption{$2C_6$ with local total antimagic 3-coloring.}\label{fig:2C6} 
\end{center}
\end{figure}
\end{example}


\begin{theorem}\label{thm-2reg} Suppose $G$ is a 2-regular graph possibly disjoint union with a path of order at least 3. If $G$ has a $C_n$ component, $n\ne 6$, then $\chi_{lt}(G)\ge 4$.  Moreover, for $m\ge 1$, $\chi_{lt}(C_3) = \chi_{lt}(mC_4)=\chi_{lt}(C_5) = 4$, and $\chi_{lt}(C_8)\le 5$. \end{theorem} 

\begin{proof} By definition, $\chi_{lt}(G)\ge 3$.  If $n\equiv 1,2\pmod{3}$, then $\chi_{lt}(G)\ge \chi_t(G)\ge 4$. Suppose $n\equiv 0\pmod{3}$ and that $\chi_{lt}(G) = \chi_t(G) = 3$. Let $f$ be a required local total antimagic 3-coloring of $G$. Let the consecutive vertices and edges of $C_n$ be $u_1,e_1,u_2,e_2,\ldots, u_n,e_n$. Suppose $n=3$. Without loss of generality, we may assume the 3 distinct weights are $a,b,c$ such that $a = w(u_{1}) = w(e_{2})$, $b=w(e_{1}) = w(u_{3})$ and $c=w(u_{2}) = w(e_{3})$. Thus, we have  

\[a=f(e_{1})+f(e_{3}) = f(u_{2})+f(u_{3})\] \[b=f(e_{2})+f(e_{3}) =f(u_{1})+f(u_{2}) \]  \[c=f(e_{1})+f(e_{2}) = f(u_{1}) + f(u_{3}).\] Now, $a+b-c$ gives $f(e_{3}) = f(u_{2})$, a contradiction. 


\ms\nt Suppose $n\ge 9$. Similarly, we may assume $a = w(u_{j}), j \equiv 1\pmod{3}, b = w(u_{j}), j\equiv 0\pmod{3}, c = w(u_{j}), j\equiv 2\pmod{3}$. Thus, we have

\[a = f(e_{3})+f(e_{4}) = f(e_{6}) + f(e_{7})\]  \[b = f(e_{2}) + f(e_{3}) = f(e_{5}) + f(e_{6})\] \[c = f(e_{1}) + f(e_{2}) = f(e_{4}) + f(e_{5}).\] From $a-b+c$, we get $f(e_{4}) + f(e_{1}) = f(e_{7}) + f(e_{4})$ so that $f(e_{1}) = f(e_{7})$, a contradiction. Thus, $\chi_{lt}(G)\ge 4$ if $G$ has a component of order $n\ne 6$.

\ms\nt Consider $C_3$, label the vertices and edges $u_{1,1}$, $e_{1,1}$, $u_{1,2}$, $e_{1,2}$, $u_{1,3}$, $e_{1,3}$ by $1,3,5,4,6,2$ bijectively so that the corresponding weights are $5,6,7,11,6,7$. Thus, $\chi_{lt}(C_3)=4$.

\ms\nt Consider $mC_4$. Since $mC_4$ has $4m$ vertices and $4m$ edges,  we define a bijection $f : V(mC_4) \cup E(mC_4) \to [1, 8m]$ such that for $1\le i\le m$,
\begin{enumerate}[(i)]
\item $f(u_{i,1}) = i$, $f(u_{i,2}) = 7m+1-i$, $f(u_{i,3}) = 3m+i$, $f(u_{i,4}) = 6m+1-i$, 
\item $f(e_{i,1}) = m+i$, $f(e_{i,2}) = 5m +1 - i$, $f(e_{i,3}) = 2m+i$, $f(e_{i,4}) = 8m+1-i$.
\end{enumerate}
Clearly, the weights of $u_{i,1},$ $e_{i,1},$ $u_{i,2},$ $e_{i,2},$ $u_{i,3},$ $e_{i,3},$ $u_{i,4},$ $e_{i,4}$ are $9m+1$, $7m+1$, $6m+1$, $10m+1$, $7m+1$, $9m+1$, $10m+1$, $6m+1$ respectively. Thus, $\chi_{lt}(mC_4)=4$.

\ms\nt For $C_5$, label the vertices and edges $u_{1,1}$, $e_{1,1}$, $u_{1,2}$, $e_{1,2}, \ldots, u_{1,5}$, $e_{1,5}$ by 1, 2, 7, 8, 5, 6, 3, 4, 9, 10 bijectively so that the corresponding weights are 12, 8, 10, 12, 14, 8, 10, 12, 14, 10. Thus, $\chi_{lt}(C_n)=4$ for $n=4,5$. For $C_8$, label the vertices and edges $u_{1,1}$, $e_{1,1}$, $u_{1,2}$, $e_{1,2}, \ldots, u_{1,8}$, $e_{1,8}$ by 1, 10, 16, 5, 2, 11, 13, 6, 3, 12, 14, 7, 4, 9, 15, 8 bijectively so that the corresponding weights are 18, 17, 15, 18, 16, 15, 19, 16, 18, 17, 19, 18, 16, 19, 17, 16. Thus, $\chi_{lt}(C_8)\le 5$.
\end{proof}


\begin{theorem}\label{thm-chilt=3}  For $m\ge 1$, $\chi_{lt}(G) = 3$ if and only if $G = P_n, n=3,6$ or $G=mC_6$  or $G=mC_6 + P_6$, $(m\ge 0)$.\end{theorem}   

\begin{proof} (Necessity) Suppose $\chi_{lt}(G)  = 3$. By Corollary~\ref{cor-2regPn}, Theorem~\ref{thm-mC6} and the proof of Theorem~\ref{thm-2reg}, we know $G=P_n, n=3,6$ or $G=mC_6$,  $G=mC_6+P_3$ or $G=mC_6 + P_6$ for $m\ge 1$. We shall prove that $\chi_{lt}(mC_6+P_3)\ne 3$.  Suppose equality holds and the $P_3$ component has vertices and edges $v_1, h_1, v_2, h_2, v_3$ consecutively.  Let $f$ be a local total antimagic 3-coloring of $mC_6+P_3$ with induced weights $a,b,c$. Without loss of generality, we may assume the weights of $u_{i,1}, e_{i,1},  \ldots, u_{i,6},$ $e_{i,6}$ of the $i$-th $C_6$ are $a,b,c$ repeatedly, and the weights of  $v_1, h_1, v_2, h_2, v_3$ are $a,b,c,a,b$ respectively.  Thus, we have \[a = f(e_{i,6})+f(e_{i,1}) = f(h_1)\]  \[b = f(e_{i,2}) + f(e_{i,3}) = f(h_2)\]  \[c = f(e_{i,1}) + f(e_{i,2}) = f(h_1) + f(h_2).\] The right hand side above gives $a+b-c=0$ so that $f(e_{i,6}) + f(e_{i,3}) = 0$, a contradiction.  Thus, $\chi_{lt}(mC_6 + P_3) \ge 4$. Therefore, $G=P_n, n=3,6$ or $G=mC_6$ or $G=mC_6+P_6$ for  $m\ge 1$. 

\ms\nt (Sufficiency) By Lemma~\ref{lem-Pn} and Theorem~\ref{thm-mC6}, we only need to show that $\chi_{lt}(mC_6 + P_6)=3$ for $m\ge 1$. Suppose the $P_6$ component has vertices and edges $v_1, h_1, \ldots, v_5, h_5, v_6$ consecutively. Since $G$ has $6m+6$ vertices and $6m+5$ edges, we define a bijection $f : V(G)\cup E(G) \to [1, 12m + 11]$ such that for $1\le i\le m$, 
\begin{enumerate}[(i)]
\item $f(v_1) = 6m +7, f(v_2) = 6m+3, f(v_3) = 6m+8, f(v_4) = 6m+4, f(v_5) = 6m+6, f(v_6) = 6m+5$,
\item $f(h_1) = 12m+11, f(h_2) = 1, f(h_3) = 12m+9, f(h_4) = 2, f(h_5) = 12m+10$,
\item $f(u_{i,1}) = 3i$, $f(u_{i,3}) = 3i+1$, $f(u_{i,5}) = 3i+2$,
\item $f(u_{i,2}) = 12m + 11 - 3i$, $f(u_{i,4}) = 12m + 9  - 3i$, $f(u_{i,6}) = 12m + 10 - 3i$,
\item $f(e_{i,1}) = 6m + 7 +  3i$, $f(e_{i,3}) = 6m+8+3i$, $f(e_{i,5}) = 6m+6+3i$,
\item $f(e_{i,2}) = 6m+3-3i$, $f(e_{i,4}) = 6m+4-3i$, $f(e_{i,6}) = 6m+5-3i$.
\end{enumerate}
Clearly, the weights of $u_{i,1},$ $e_{i,1}, \ldots, u_{i,6},$ $e_{i,6}$ are $12m+12$, $12m+11$, $12m+10$ repeatedly for $1\le i\le m$ whereas the weights of $v_1, h_1, \ldots, v_5, h_5, v_6$ are $12m+11, 12m+10, 12m+12, 12m+11, 12m+10, 12m+12, 12m+11, 12m+10$ respectively. Thus, $f$ is a local total antimagic labeling and $\chi_{lt}(G) = 3$. 
\end{proof}

\begin{example} The figure below gives the local total antimagic 3-coloring of $2C_6+P_6$ as defined above with induced weights $34,35,36$.
\begin{figure}[H]
\begin{center}
\centerline{\epsfig{file=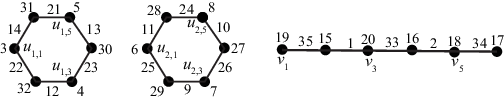, width=9.5cm}}
\caption{$2C_6+P_6$ with local total antimagic 3-coloring.}\label{fig:2C6P6} 
\end{center}
\end{figure}
\end{example}


\begin{corollary} Suppose $m\ge 1$ is odd and $s\ge \frac{m+1}{2}$. If $G_v(s)$ is obtained from $mC_4$ by attaching $s\ge 2$ pendant edges to a single vetex of $mC_4$, then $\D(G_v(s))+1=s+3\le \chi_{lt}(G_v(s)) \le s+4=\D(G_v(s))+2$. \end{corollary}

\begin{proof} Let $f$ be the local total antimagic labeling of $mC_4$ as defined in the proof of Theorem~\ref{thm-2reg}. Without loss of generality, assume the $s$ edges are attached to $u_{1,1}$. Define a labeling $g$ of $G_v(s)$ such that $g(x) = f(x)$ if $x\in V(G)\cup E(G)$, and label the vertices and edges of the $s$ added pendant edges as in the proof of Theorem~\ref{thm-HfromG} for $k=1$. Now, $g$ is a local total antimagic labeling that induces weights $9m+1,7m+1,6m+1,10m+1$,  $8m+2i, 1\le i\le s$, $\sum^s_{i=1} (8m+2i) = 8ms+s(s+1)$. Since $s\ge \frac{m+1}{2}$, $9m+1\in\{8m+2i\,|\, 1\le i\le s\}$, there are exactly $s+4$ distinct weight so that $\chi_{lt}(G_v(s))\le s+4$. Since $\D(G_v(s)) = s+2$ and $G_v(s)$ has $s$ pendant edges, by Corollary~\ref{cor-pendant}, $\chi_{lt}(G_v(s))\ge s+3$. This completes the proof.
\end{proof}

\section{$\chi_{lt} = $ number of pendant edges + 1}\label{sec-pend+1}

\nt We first note that applying Theorem~\ref{thm-HfromG} to the vertex of $mC_6, m\ge 1$ with label 1 as in Theorem~\ref{thm-mC6}, we get the following corollary.

\begin{corollary}\label{cor-G+pend} Suppose $m\ge 1$. If $G_v(s)$ is obtained from $mC_6$ by attaching $s\ge 2$ pendant edges to a single vertex of $mC_6$, then $\chi_{lt}(G_v(s)) = \D(G_v(s))+1 = s+3$. \end{corollary}

\nt In\cite{S+N+L}, the authors proved that $\chi_{lt}(nP_3) =2n+1$ for $n\ge 1$. We can now extend the obtained labeling to the following theorem. \\


\begin{theorem}\label{thm-mC6nP3} For $m\ge 1, n\ge 2$, $\chi_{lt}(mC_6 + P_3) = 4$ and $\chi_{lt}(mC_6+nP_3)= 2n+1$.
\end{theorem}

\begin{proof}  Consider $G=mC_6 + P_3$. From the proof of Theorem~\ref{thm-chilt=3}, we know that $\chi_{lt}(G)\ge 4$.  
%
Since $G$ has $6m+3$ vertices and $6m+2$ edges, we define a bijection $f:V(G)\cup E(G)\to[1,12m+5]$ such that for $1\le i\le m$,
\begin{enumerate}[(i)]
\item $f(v_1)=1, f(h_1)=12m+5, f(v_2) = 12m+3, f(h_2)=12m+4, f(v_3)=2$,
\item $f(u_{i,1}) = 3i, f(u_{i,3}) = 3i+1, f(u_{i,5}) = 3i+2$,
\item $f(u_{i,2}) = 12m+5-3i, f(u_{i,4}) = 12m+3-3i, f(u_{i,6}) = 12m+4-3i$,
\item $f(e_{i,1}) = 6m+1+3i, f(e_{i,3}) = 6m+2+3i, f(e_{i,5}) = 6m+3i$,
\item $f(e_{i,2}) = 6m+3-3i, f(e_{i,4}) = 6m+4-3i, f(e_{i,6}) = 6m+5-3i$.
\end{enumerate}
Clearly, the weights of  $u_{i,1},$ $e_{i,1},$ $u_{i,2},$ $e_{i,2},$ $u_{i,3},$ $e_{i,3},$ $u_{i,4},$ $e_{i,4},$ $u_{i,5},$ $e_{i,5},$ $u_{i,6},$ $e_{i,6}$ are 
$12m+6$, $12m+5$, $12m+4$ repeatedly for $1\le i\le m$ whereas the weights of $v_1, h_1, v_2, h_2, v_3$ are $12m+5, 12m+4, 24m+9, 12m+5, 12m+4$ respectively. Thus, $\chi_{lt}(G) \le 4$. Consequently, $\chi_{lt}(G) = 4$.

\ms\nt Consider $n\ge 2$. Now, $G=mC_6+nP_3$ has $2n+1$ pendant vertices and maximum degree 3. By Corollary~\ref{cor-pendant}, we have $\chi_{lt}(G)\ge 2n+1$. Suppose the $nP_3$ has vertex set $\{v_{j,1}, v_{j,2}, v_{j,3}\}$ and edge set $\{v_{j,1}v_{j,2}, v_{j,2}v_{j,3}\}$ for $1\le j\le n$.  Define a total labeling $f : V(G) \cup E(G) \to [1, 12m+5n]$  such that for $1\le i\le m, 1\le j\le n$,
\begin{enumerate}[(a)]
\item $f(v_{j,1})=j, f(v_{j,1}v_{j,2}) = 5n+12m+1-j, f(v_{j,2}) = 3n+12m+1-j, f(v_{j,2}v_{j,3}) = 3n+12m+j, f(v_{j,3}) = n+j$,
\item $f(u_{i,1}) = 2n+3i - 2$, $f(u_{i,3}) = 2n+3i-1$, $f(u_{i,5}) = 2n+3i$,
\item $f(u_{i,2}) = 2n+12m + 3 - 3i$, $f(u_{i,4}) = 2n+12m + 1 - 3i$, $f(u_{i,6}) = 2n+12m + 2 - 3i$,
\item $f(e_{i,1}) = 2n+6m - 1 + 3i$, $f(e_{i,3}) = 2n+6m + 3i$, $f(e_{i,5}) = 2n+6m - 2 + 3i$,
\item $f(e_{i,2}) = 2n+6m +1 - 3i$, $f(e_{i,4}) = 2n+6m+2-3i$, $f(e_{i,6}) = 2n+6m+3-3i$. 
\end{enumerate}
We now have $w(v_{j,1}) = 5n+12m+1-j$, $w(v_{j,3})=3n+12m+j$ for $1\le j\le n$, $w(v_{j,2}) = 8n+24m+1$, $w(v_{j,1}v_{j,2}) = 3n+12m+1, w(v_{j,2}v_{j,3}) = 4n+12m+1$ so that the $nP_3$ components have $2n+1$ distinct weights. Moreover, the weights of $u_{i,1}$, $e_{i,1}$, $u_{i,2}$, $e_{i,2}$, $u_{i,3}$, $e_{i,3}$, $u_{i,4}$, $e_{i,4}$, $u_{i,5}$, $e_{i,5}$, $u_{i,6}$, $e_{i,6}$ are $4n+12m+2, 4n+12m+1, 4n+12m$ (that are the $w(v_{n-1,1}), w(v_{n,1}), w(v_{n,3})$) repeatedly for $1\le i\le m$.  Thus, $f$ is a local total antimagic labeling that induces $2n+1$ distinct weights so that $\chi_{lt}(G) \le 2n+1$. Consequently, $\chi_{lt}(G)  = 2n+1$. 
\end{proof}

\begin{example} The figure below gives the local total antimagic 5-coloring of $2C_6+2P_3$ with induced weights $31,32,33,34,65$ as defined above.
\begin{figure}[H]
\begin{center}
\centerline{\epsfig{file=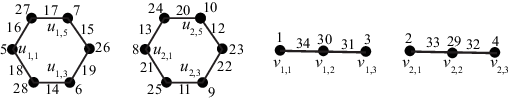, width=9.5cm}}
\caption{$2C_6+2P_3$ with local total antimagic 5-coloring.}\label{fig:2C62P3} 
\end{center}
\end{figure}
\end{example}

\nt In~\cite{S+N}, the authors also proved that $\chi_{lt}(mP_6) = 2m+1$ for $m\ge 2$, and that $\chi_{lt}(mP_6 + nP_3) = 2m+2n+1$ for $m\ge 1, n\ge 2$. We can now prove the following theorems.

\begin{theorem}\label{thm-mC6+nP6} For $m, n\ge 1$, $\chi_{lt}(mC_6+nP_6) = 2n+1$.   \end{theorem}  

\begin{proof} Let $G=mC_6+nP_6$. By Theorem~\ref{thm-chilt=3}, we only need to consider $m\ge 1, n\ge 2$. Since $G$ has $2n+1$ pendant edges and $\D(G)=2$, by Corollary~\ref{cor-pendant}, $\chi_{lt}(G)\ge 2n+1$. We shall show that $\chi_{lt}(G)\le 2n+1$. Keep the notations for the $mC_6$, $m\ge 1$. For $1\le t\le n$, let the vertex and edge sets of the $t$-th component of the $nP_6$ be $\{y_{t,a}\mid 1\le a\le 6\}$ and $\{z_{t,a} = y_{t,a}y_{t,a+1}\mid 1\le a\le 5\}$ respectively. Since $G$ has order $6m+6n$ and size $6m+5n$, we define a bijection $f : V(G) \cup E(G) \to [1, 12m+11n]$ such that for $1\le i\le m$, $1\le t\le n$,
\begin{enumerate}[(a)]
\item $f(u_{i,1}) = 2n + 3i - 2$, $f(u_{i,3}) = 2n + 3i - 1$, $f(u_{i,5}) = 2n + 3i$,
\item $f(u_{i,2}) = 8n + 12m - 3i + 3$, $f(u_{i,4}) = 8n + 12m - 3i + 1$, $f(u_{i,6}) = 8n + 12m - 3i + 2$,
\item $f(e_{i,1}) = 8n + 6m + 3i - 1$, $f(e_{i,3}) = 8n + 6m+3i$, $f(e_{i,5}) = 8n + 6m + 3i + 1$,
\item $f(e_{i,2}) = 2n+6m - 3i + 1$, $f(e_{i,4}) = 2n + 6m - 3i + 2$, $f(e_{i,6}) = 2n + 6m - 3i + 3$,
\item $f(y_{t,1})  = 8n + 6m - 3t + 2$, $f(y_{t,3}) = 8n + 6m - 3t + 3$, $f(y_{t,5}) = 8n + 6m - 3t + 1$,
\item $f(y_{t,2}) = 2n + 6m + 3t - 2$, $f(y_{t,4}) = 2n + 6m + 3t -1$, $f(y_{t,6}) = 2n+6m+3t$,
\item $f(z_{t,1}) = 11n + 12m + 1 - t$, $f(z_{t,3}) = 9n + 12m + 1- t$, $f(z_{t,5}) = 10n + 12m + 1 - t$,
\item $f(z_{t,2}) = t$, $f(z_{t,4}) = n+t$.
\end{enumerate}
One can check that the induced weights of each component of the $mC_6$, starting from $w(u_{i,1})$ and $w(e_{i,1})$, are $10n+12m+2, 10n+12m+1, 10n+12m$ repeatedly and the induced non-pendant vertex weights of each component of the $nP_6$, starting from $w(z_{t,1})$, are $10n+12m$, $11n+12m+1$, $10n+12m+1$, $9n+12m+1$, $10n+12m+2$, $10n+12m+1$, $10n+12m$, $11n+12m+1$, $10n+12m+1$ consecutively, while $w(y_{t,1}) = 11n+12m+1-t$ and $w(y_{t,6}) = 10n+12m+1-t$ for $1\le t\le n$. Thus, $f$ is a local total antimagic labeling that induces $2n+1$ distinct weights so that $\chi_{lt}(G)\le 2n+1$. This completes the proof.
\end{proof}

\begin{example} The figure below gives the local total antimagic 5-coloring of $2C_6+2P_6$ as defined above with induced weights $43,44,45,46,47$.
\begin{figure}[H]
\begin{center}
\centerline{\epsfig{file=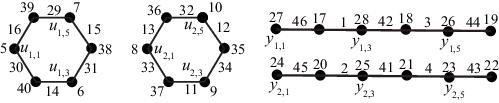, width=9.5cm}}
\caption{$2C_6+2P_6$ with local total antimagic 5-coloring.}\label{fig:2C62P6} 
\end{center}
\end{figure}
\end{example}

\begin{theorem}\label{thm-mC6+nP6+aP3} For $m,n\ge 1, a\ge 2$,if $n\ge 2a+2$ or else $a\ge 2n$, then $\chi_{lt}(mC_6+nP_6+aP_3)= 2n+2a+1$. Otherwise, 
\begin{enumerate}[(i)]
\item if $n=2a+1$ or else $s\in[n+1,2n-1]$, then  $2n+2a+1\le \chi_{lt}(mC_6+nP_6+aP_3)\le 2n+2a+2$. 
\item if $n=2a$ or else $a = n$, then  $2n+2a+1\le \chi_{lt}(mC_6+nP_6+aP_3)\le 2n+2a+3$. 
\item if $n\in [a+1,2a-1]$ or else $a=n-1$, then  $2n+2a+1\le \chi_{lt}(mC_6+nP_6+aP_3)\le 2n+2a+4$. 
\item if $n\le a$ or else $a\le n-2$, then  $2n+2a+1\le \chi_{lt}(mC_6+nP_6+aP_3)\le 2n+2a+5$. 
\end{enumerate}
\end{theorem}

\begin{proof} Let $G = mC_6 + nP_6 + aP_3$ for  $m,n\ge 1, a\ge 2$. Since $G$ has $2n+2a$ pendant edges and $\D(G) = 2$, by Corollary~\ref{cor-pendant}, $\chi_{lt}(G)\ge 2n+2a+1$. We shall show that $\chi_{lt}(G)$ has the given upper bound under the given conditions. Keep the notations of $mC_6 + nP_6$ as defined in the proof of Theorem~\ref{thm-mC6+nP6}. Let $f$ be the total labeling as defined in the proof of Theorem~\ref{thm-mC6+nP6}. For $1\le r\le a$, let the vertex and edge sets of the $r$-th component of the $aP_3$ be $\{v_{r,b}\mid1\le b\le 3\}$ and $\{h_{r,b} = v_{r,b}v_{r,b+1}\mid 1\le b\le 2\}$ respectively.  

\ms\nt Since $G$ has order $6m+6n+3a$ and size $6m+5n+2a$, we define a bijection $g:V(G) \cup E(G) \to [1,12m+11n+5a]$ such that for $1\le i\le m, 1\le t\le n, 1\le r\le a$, $g(x) = f(x)+2a$ for $x\in V(mC_6+nP_6)\cup E(mC_6+nP_6)$. Moreover, $g(v_{r,1}) = r, g(h_{r,1}) = 12m+11n+5a+1-r$, $g(v_{r,2}) = 12m+11n+3a+1-r$, $g(h_{r,2}) = 12m+11n+3a+r$ and $g(v_{r,3}) = a+r$.  

\ms\nt Similar to the proof of Theorem~\ref{thm-mC6+nP6}, one can check that the induced weights of each component of the $mC_6$ are $10n+12m+4a+2$, $10n+12m+4a+1$, $10n+12m+4a$ repeatedly, the induced non-pendant vertex weights of each component of the $nP_6$ are $10n+12m+4a$, $11n+12m+4a+1$, $10n+12m+4a+1$, $9n+12m+4a+1$, $10n+12m+4a+2$, $10n+12m+4a+1$, $10n+12m+4a$, $11n+12m+4a+1$, $10n+12m+4a+1$, while $w_g(y_{t,1}) = 11n+12m+2a+1-t$, $w_g(y_{t,5}) = 10n+12m+2a+1-t$, $w_g(v_{r,1}) = 11n+12m+5a+1-r$, $w_g(h_{r,1}) = 11n+12m+3a+1$, $w_g(v_{r,2}) = 22n+24m+8a+1$, $w_g(h_{r,2})=11n+12m+4a+1$ and $w_g(v_{r,3}) = 11n+12m+3a+r$ for $1\le t\le n, 1\le r\le a$. Clearly, $w_g(v_{r,2})=22n+24m+8a+1$ is not in the pendant vertex weights set $W=[9n+12m+2a+1, 11n+12m+2a] \cup [11n+12m+3a+1,11n+12m+5a]$. Note that this labeling is local total antimagic. Moreover, there are exactly $2n+2a+1$ distinct weights if and only if all the induced weights set of the non-pendant vertices of $mC_6+nP_6$ components, namely $\{9n+12m+4a+1,10n+12m+4a,10n+12m+4a+1, 10n+12m+4a+2, 11n+12m+4a+1\}$ is a subset of $W$. This implies that 

\begin{enumerate}[(a)]
\item $2n\ge 2a+1$, or else $a\ge 2n$ for $9n+12m+4a+1 \in W$,
\item $n\ge 2a$, or else $a\ge n+1$ for $10n+12m+4a\in W$,
\item $n\ge  2a+1$, or else $a\ge n$ for $10n+12m+4a+1\in W$,\
\item $n \ge 2a+2$ , or else $a\ge n-1$ for $10n+12m+4a+2\in W$.   
\end{enumerate}
Therefore, if $n\ge 2a+2$ or else $a\ge 2n$, then $\chi_{lt}(mC_6+nP_6+aP_3)\le 2n+2a+1$. Otherwise, 
\begin{enumerate}[(i)]
\item if $n=2a+1$ or else $a\in[n+1,2n-1]$, then there are $2n+2a+2$ distinct weights so that $2n+2a+1\le \chi_{lt}(mC_6+nP_6+aP_3)\le 2n+2a+2$. 
\item if $n=2a$ or else $a = n$, then there are $2n+2a+3$ distinct weights so that $2n+2a+1\le \chi_{lt}(mC_6+nP_6+sP_3)\le 2n+2a+3$. 
\item if $n\in [a+1,2a-1]$ or else $a=n-1$, then there are $2n+2a+4$ distinct weights so that $2n+2a+1\le \chi_{lt}(mC_6+nP_6+aP_3)\le 2n+2s+4$. 
\item if $n\le a$ or else $a\le n-2$, then there are $2n+2a+5$ distinct weighst so that $2n+2a+1\le \chi_{lt}(mC_6+nP_6+aP_3)\le 2n+2a+5$. 
\end{enumerate}
This completes the proof. 
\end{proof}

\begin{example} The figure below gives the local total antimagic 11-coloring of $C_6+4P_6+P_3$ with $n=4,a=1$ and induced weights $51, 52, \ldots, 58, 60, 61, 121$, as well as local total antimagic 9-coloring of $C_6+P_6+3P_3$ with $n=1,a=3$ and induced weights $28,29,33,34,\ldots,38, 71$ as defined above.
\begin{figure}[H]
\begin{center}
\centerline{\epsfig{file=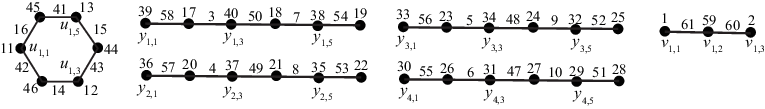, width=13.3cm}}
\caption{$C_6+4P_6+P_3$ with local total antimatic 13-coloring.}\label{fig:C64P62P3} 
\vskip1cm
\centerline{\epsfig{file=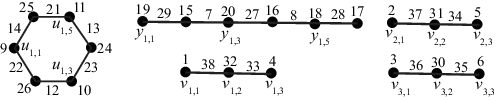, width=9.5cm}}
\caption{$C_6+P_6+3P_3$ with local total antimagicl 9-coloring.}\label{fig:C6P63P3} 
\end{center}
\end{figure}
\end{example}

\begin{theorem}\label{thm-mC6+nP3+ks} Let $G = mC_6 + nP_3$, $m\ge 1, n\ge 2$. Suppose $s\ge 1, 1\le j'\le s, 1\le i'\le k, 1\le i\le m, 1\le j\le n$, $k\in\{2n+1\}\cup [2n+12m+1,3n+12m]$, then $ks+2n+1 \le \chi_{lt}(G_v(k,s)) \le ks+2n+2$. Moreover, $\chi_{lt}(G_v(k,s)) = ks+2n+1$ if
\begin{enumerate}[(a)]
\item  $v=v_{1,1}$, $k=1$, $n$ is odd and $s\ge (3n+24m+1)/2$,
\item $v=v_{j,1}$ for $j\in[2,n]$, $k=j$ and $2(j'-1)j+1 \le 3n+12m+1\le 2j'j$,
\item $v=v_{j,3}$ for $j\in [1,n]$, $k=n+j$ and $(2j'-1)(n+j)+1 \le  3n+12m+1 \le 2j'(n+j)$,
\item  $v=v_{j,2}$, $k=3n+12m+1-j, s\ge 1$,
\item  $v=u_{i,1}$, $k=2n+3i-2$, and $8n+24m+1\in\bigcup^s_{j'=1}\,[5n+12m+(2j'-1)(2n+3i-2)+1, 5n+12m+2j'(2n+3i-2)]$,
\item  $v=u_{i,3}$, $k=2n+3i-1$, and  $8n+24m+1\in\bigcup^s_{j'=1}\, [5n+12m+(2j'-1)(2n+3i-1)+1, 5n+12m+2j'(2n+3i-1)]$,
\item  $v=u_{i,5}$, $k=2n+3i$, and $8n+24m+1\in\bigcup^s_{j'=1}\,[5n+12m+(2j'-1)(2n+3i)+1,5n+12m+2j'(2n+3i)]$,
\item  $v=u_{i,2}$, $k=2n+12m+3-3i, s\ge 1$, 
\item  $v=u_{i,4}$, $k=2n+12m+1-3i, s\ge 1$, 
\item  $v=u_{i,6}$, $k=2n+12m+2-3i, s\ge 1$. 
\end{enumerate}
\end{theorem}

\begin{proof} Let $f$ be the local total antimagic labeling of $G$ as defined in the proof of Theorem~\ref{thm-mC6nP3}. Clearly, if $v$ is a pendant vertex, then $v\in \{v_{j,1}, v_{j,3}\mid 1\le j\le n\}$ and $f(v) = k = j$ for $v=v_{j,1}$ and $k=n+j$ for $v=v_{j,3}$. Moreover, if $v$ is not a pendant vertex, then $v\in\{u_{i,1}, \ldots, u_{i,6}, v_{j,2}\,|\, 1\le i\le m, 1\le j\le n\}$ and $f(v) = k \in\{2n+3i-2, 2n+3i-1, 2n+3i, 2n+12m+3-3i, 2n+12m+2-3i, 2n+12m+1-3i, 3n+12m+1-j\,|\,1\le i\le m, 1\le j\le n\}$. Moreover, $\D(G_v(k,s)) = ks+2$ and $G_v(k,s)$ has $ks+2n$ pendant edges. Suppose the added $ks$ pendant edges incident to $v$ are $e_{i',j'}$ and the corresonding pendant vertices are $x_{i',j'}$ for $1\le j'\le s, 1\le i'\le k$. By Corollary~\ref{cor-pendant}, $\chi_{lt}(G_v(k,s))\ge ks+2n+1$. 

\ms\nt Since $G_v(k,s)$ has $3n+6m+ks$ vertices and $2n+6m+ks$ edges, we define a total labeling $g: V(G_v(k,s)) \cup E(G_v(k,s)) \to [1, 5n+12m+2ks]$ such that $g(z) = f(z)$ for $z\in V(G)\cup E(G)$. Otherwise, $g(z)$ is as defined in the proof of Theorems~\ref{thm-HfromG} or~\ref{thm-HfromG2} . One can check that $w_g(v) = w_f(v) + \sum^s_{j'=1} \sum^k_{i'=1} (12m+5n + (2j'-1)k+i')$, and $w(z)\in [3n+12m+1,5n+12m]\cup\{8n+24m+1\}\cup \{5n+12m+(2j'-1)k+i'\,|\,1\le j'\le s, 1\le i'\le k\}$ for $z\in (V(G_v(k,s)) \cup E(G_v(k,s))) \setminus \{v\}$. Moreover, $g$ is a local total antimagic labeling. Thus, $\chi_{lt}(G_v(k,s))\le ks+2n+2$ if all the weights are distinct. We shall need to check the conditions under Theorem~\ref{thm-HfromG2} in the following 10 cases.

\begin{enumerate}[(a)]
\item Suppose $v=v_{1,1}$ so that $k=1$. Thus, $\{w_g(e_{j',1})\} = \{5n+12m+2j' \mid 1\le j' \le s, s\ge 2\}$.  Therefore, $8n+24m+1 \in \{w_g(e_{j',1})\}$ so that $n$ is odd and $s\ge (3n+24m+1)/2$. 

\item Suppose $v=v_{j,1}, 2\le j\le n$ so that $k = j\ge 2$. Therefore, $8n+24m+1 \in \{w_g(e_{j',i'})\} = \bigcup^s_{j'=1} [5n+12m + (2j'-1)j+1, 5n+12m+2j'j]$, $2\le j\le n$, so that $2(j'-1)j+1 \le 3n+12m+1\le 2j'j$. As an example, take $n=3,m=1$, we can choose $j=3, s\ge j'=4$ to get $8n+24m+1=49\in [49, 52]=[5n+12m + (2j'-1)j+1, 5n+12m+2j'j]$ as required.

\item Suppose $v=v_{j,3}, 1\le j\le n$ so that $k = n+j$.  Therefore, $8n+24m+1 \in \{w_g(e_{j',i'})\} = \bigcup^s_{j'=1} [5n+12m + (2j'-1)(n+j)+1, 5n+12m+2j'(n+j)]$, $1\le j\le n$, so that  $(2j'-1)(n+j)+1 \le  3n+12m+1 \le 2j'(n+j)$. As an example, take $n=m=3$, we can choose $j=1, s\ge j'=6$ to get $8n+24m+1=97\in [96, 99] =  [5n+12m + (2j'-1)(n+j)+1, 5n+12m+2j'(n+j)]$ as required.

\item Suppose $v=v_{j,2}$ so that $k=3n+12m+1-j$. Thus, $\{w_g(e_{j',i'})\} = \{5n+12m+(2j'-1)(3n+12m+1-j)+i' \mid 1\le i'\le k, 1\le j'\le s\} = \bigcup^s_{j'=1}\, [5n+12m+(2j'-1)(3n+12m+1-j)+1, 5n+12m+2j'(3n+12m+1-j)]$. If $j'=1$, then $8n+24m+1\in [8n+24m+2-j, 11n+36m+2-2j]$ since $j\le n$. Thus, $g$ induces $ks+2n+1$ distinct weights. Therefore, if $G=mC_6+nP_3$, then $\chi_{lt}(G_{v_{j,2}}(3n+12m+1-j,s)) =  (3n+12m+1-j)s+2n+1$ for $s\ge 1$.



\item Suppose $v=u_{i,1}$ so sthat $k=2n+3i-2$. Thus, $\{w_g(e_{j',i'})\} = \{5n+12m+(2j'-1)k+i'\,|\, 1\le i' \le k, 1\le j'\le s\} = \bigcup^s_{j'=1}\, [5n+12m+(2j'-1)(2n+3i-2)+1, 5n+12m+2j'(2n+3i-2)]$, denoted $U$, for $1\le i\le m, 1\le j'\le s$ Thus, if $8n+24m+1\in U$, then $g$  induces $ks+2n+1$ distinct weights.   As an example: take $m=1, n=3, i=1, k=7, j'=2$, we get $8n+24m+1=49$ and $[5n+12m+(2j'-1)(2n+3i-2)+1, 5n+12m+2j'(2n+3i-2)] = [49,55]$. Thus, for $G=C_6 + 3P_3$,  $\chi_{lt}(G_{u_{1,1}}(7,s)) = 7s+7$ for $s\ge 2$. 

\item Suppose $v=u_{i,3}$ so that $k=2n+3i-1$. By a similar argument, we get $g$ induces $ks+2n+1$ distinct weights if $8n+24m+1\in  \bigcup^s_{j'=1}\, [5n+12m+(2j'-1)(2n+3i-1)+1, 5n+12m+2j'(2n+3i-1)]$ for $1\le i\le m, 1\le j'\le s$. As an example: take $m=2, n = 3, i=2, k=11, j'=2$, we get $8n+24m+1=73$ and $[5n+12m+(2j'-1)(2n+3i-1)+1, 5n+12m+2j'(2n+3i-1)] = [73,83]$. Thus, for $G=2C_6 + 3P_3$, $\chi_{lt}(G_{u_{2,3}}(11,s))=11s+7$ for $s\ge 2$.

\item Suppose $v=u_{i,5}$ so that $k=2n+3i$. By a similar argument, we get $g$ induces $ks+2n+1$ distinct weights if $8n+24m+1\in \bigcup^s_{j'=1}\,[5n+12m+(2j'-1)(2n+3i)+1,5n+12m+2j'(2n+3i)]$. As an example: take $m= n= i = 3, k=15, j' = 2$, we get $8n+24m+1 = 97$ and $[5n+12m+(2j'-1)(2n+3i)+1,5n+12m+2j'(2n+3i)] = [97, 111]$. Thus, for $G=3C_6 + 3P_3$, $\chi_{lt}(G_{u_{3,5}}(15,s)) = 15s+7$ for $s\ge 2$.

\item Suppose $v=u_{i,2}$ so that $k=2n+12m+3-3i$. We get $g$ induces $ks+2n+1$ distinct weights if $8n+24m+1\in\bigcup^s_{j'=1}\,[5n+12m+(2j'-1)(2n+12m+3-3i) +1, 5n+12m+2j'(2n+12m+3-3i)]$.  For $j'=1$,  we get $8n+24m+1 \in [7n+12m+4-3i, 9n+36m+6-6i]$ since $i\le m$. Thus, for $G=mC_6 + nP_3$, $\chi_{lt}(G_{u_{i,2}}(2n+12m+3-3i,s)) = (2n+12m+3-3i)s+2n+1$ for $s\ge 1$.

\item Suppose $v=u_{i,4}$ so that $k=2n+12m+1-3i$. We get $g$ induces $ks+2n+1$ distinct weights if $8n+24m+1\in\bigcup^s_{j'=1}\,[5n+12m+(2j'-1)(2n+12m+1-3i) +1, 5n+12m+2j'(2n+12m+1-3i)]$. For $j'=1$, we get $8n+24m+1\in [7n+24m+2-3i, 9n+36m+3-6i]$ since $i\le m$. Thus, for $G=mC_6+nP_3$,  $\chi_{lt}(G_{u_{i,4}}(2n+12m+1-3i,s)) = (2n+12m+1-3i)s + 2n+1$ for $s\ge 1$. 

\item Suppose $v=u_{i,6}$ so that $k=2n+12m+2-3i$. We get $g$ induces $ks+2n+1$ distinct weights if $8n+24m+1\in\bigcup^s_{j'=1}\,[5n+12m+(2j'-1)(2n+12m+2-3i) +1, 5n+12m+2j'(2n+12m+2-3i)]$. For $j'=1$, we get $8n+24m+1\in [7n+24m+3-3i, 9n+36m+4-6i]$ since $i\le m$. Thus, for $G=mC_6 + nP_3$, $\chi_{lt}(G_{u_{i,6}}(2n+12m+2-3i,s)) = (2n+12m+2-3i)s + 2n+1$ for $s\ge 1$. 
\end{enumerate}

\ms\nt This completes the proof.   \end{proof}

\begin{example} The figure below gives the local total antimagic 65-coloring of $G_v(k,s)$ for $G=2C_6+2P_3, v = v_{1,2}, k=30, s=2$ with induced weights $31,32,33,34,65,66,\ldots, 94,$ $125,126,\ldots,$ $154$ and  $6635$ being the weight of $v_{1,2}$ in $G_v(k,s)$ as defined above.

\begin{figure}[H]
\begin{center}
\centerline{\epsfig{file=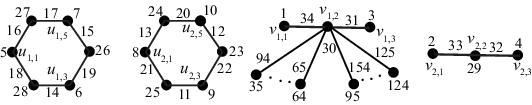, width=12cm}}
\caption{$G_v(k,s)$, $v=v_{1,2}$, $k=30,s=2$ with local total antimatic 65-coloring.}\label{fig:2C62P3+ks} 
\end{center}
\end{figure}
\end{example}

\nt Note that Theorem~\ref{thm-HfromG2} is not applicable to $mC_6 + nP_6$ for $m, n\ge 1$. Let $G=mC_6 + nP_6 + aP_3$ with vertex and edge sets as defined in the proof of Theorem~\ref{thm-mC6+nP6+aP3}. We now have the following theorem.

\begin{theorem}\label{mC6+nP6+aP3+ks} For $m,n\ge 1, a\ge 2$ and $G=mC_6 + nP_6 + aP_3$,  if $n\ge 2a+2$ or else $a\ge 2n$, then $\chi_{lt}(G_v(k,s)) = 2n+2a+ks+1$ for 
\begin{enumerate}[(1)]
\item $v=y_{t,1}$, $k=8n+6m-3t+2a+2$, $t\le (5n+3)/6$,
\item $v=y_{t,6}$, $k=2n+6m+3t+2a$ and 
    \begin{enumerate}[(a)]
	\item $s\ge 1, t\in[1,n], a\ge 2n$ only, or else 
	\item $s\ge 2, (3n-12m-5a+1)/12\le t\le (5n-6m-3a+1)/9$, or else
	\item $s\ge 3, t \le (n-18m-7a+1)/15$.
   \end{enumerate}
\item $v=v_{r,1}, r\in [1,a], k=r, (2j'-1)r \le 11n+12m+3a+1\le 2j'r$,
\item $v=v_{r,3}, r\in [1,a], k=a+r, (2j'-1)(a+r) \le 11n+12m+3a+1\le 2j'(a+r)$,
\item $v=u_{i,1}, i\in [1,m], k=2n+2a+3i-2, (2j'-1)(2n+2a+3i-2) \le 11n+12m+3a+1 \le 2j'(2n+2a+3i-2)$,
\item $v=u_{i,3}, i\in [1,m], k=2n+2a+3i-1, (2j'-1)(2n+2a+3i-1) \le 11n+12m+3a+1 \le 2j'(2n+2a+3i-1)$,
\item $v=u_{i,5}, i\in [1,m], k=2n+2a+3i, (2j'-1)(2n+2a+3i) \le 11n+12m+3a+1 \le 2j'(2n+2a+3i)$,
\item $v=u_{i,2}, i\in [1,m], k = 8n+12m+2a-3i+3, s\ge 1$,
\item $v=u_{i,4}, i\in [1,m], k = 8n+12m+2a-3i+1, s\ge 1$, 
\item $v=u_{i,6}, i\in [1,m], k = 8n+12m+2a-3i+2, s\ge 1$,  
\item $v=y_{t,3},k = 8n+6m+2a-3t+3, t\le (5n+a+5)/6$, $n\ge a+5$ for $a=1,2$ and $n\ge 2a+2$ otherwise,
\item $v=y_{t,5},k = 8n+6m+2a-3t+1, t\le (5n+a+1)/6$ and $n\ge 2a+2$ only,
\item $v=y_{t,2},k = 2n+6m+2a+3t-2$ and 
     \begin{enumerate}[(a)]
	\item $s\ge 1, t\in [1,n], a\ge n+5$ for $a=1,2,3,4$ and $a\ge 2n$ otherwise, or else
	\item $s\ge 2, (3n-12m-5a+9)/12 \le t\le (5n-6m-a+7)/3$, or else 
	\item $s\ge 3, t\le (n-18m-7a+13)/15$,
    \end{enumerate}
\item $v=y_{t,4}$, $k=2n+6m+2a+3t-1$, and
    \begin{enumerate}[(a)]
	\item $s\ge 1, t\in[1,n], a\ge n+3$ for $a=1,2$ and $a\ge 2n$ otherwise, or else
	\item $s\ge 2, (3n-12m-5+5)/12 \le t \le (5n-6m-3a+4)/9$, or else
	\item $s\ge 3, t\le (n-18m-7a+6)/15$.
    \end{enumerate}
\item $v=v_{r,2}, r\in[1,a]$, $k=12m+11n+3a+1-r, s\ge 1$. 
\end{enumerate}
\end{theorem}

\begin{proof} Let $g$ be the local total antimagic labeling of $G$ as defined in the proof of Theorem~\ref{thm-mC6+nP6+aP3}. By an argument similar to the proof of Theorem~\ref{thm-mC6+nP3+ks} and Theorem~\ref{thm-HfromG2}(iii), we shall show that the only weight $w_g(v_{r,2})=22n+24m+8a+1$ which is not a pendant edge label under $g$ must be in $\{w_g(e_{j',i'})\} = \bigcup^s_{j'=1} [11n+12m+5a+(2j'-1)k+1, 11n+12m+5a+2j'k]$. We also check the conditions under Theorem~\ref{thm-HfromG2}. We have the following first four cases for $v$ being a pendant vertex and next 11 cases for $v$ not a pendant vertex.
\begin{enumerate}[(1)]
\item Suppose $v=y_{t,1}, t\in[1,n]$ with $k=8n+6m-3t+2a+2$. Therefore, $22n+24m+8a+1 \in \{w_g(e_{j',i'})\}$ means $(2j'-1)(8n+6m-3t+2a+2)+1 \le 11n+12m+3a+1 \le 2j'(8n+6m-3t+2a+2)$. Thus, $j'=1$ so that $8n + 6m - 3t + 2a+3 \le 11n+12m+3a+1 \le 16n+12m-6t+4a+4$. Consequently, $t\le (5n+3)/6$.
\item Suppose $v=y_{t,6}, t\in[1,n]$ with $k=2n+6m+3t+2a$. Therefore, we must have $(2j'-1)(2n+6m+3t+2a) \le 11n+12m+3a+1\le 2j'(2n+6m+3t+2a)$. Thus, $j'\le 3$. \\ If $j'=1$, then $2n+6m+3t+2a \le 11n+12m+3a+1\le 4n+12m+6t+4a$ so that $(7n-a+1)/6 \le t\le (9n+6m+a+1)/3$. Since $t\le n$, we then have $a\ge n+1$. The conditions $n\ge 2a+2$ or else $a\ge 2n$ further implies that we must have $a\ge 2n$ for this case. \\ If $s\ge j'=2$, then $6n+18m+9t+6a\le 11n+12m+3a+1\le 8n + 24m+12t+8a$ so that $(3n-12m-5a+1)/12\le t\le (5n-6m-3a+1)/9$.  \\ If $s\ge j'=3$, then $10n+30m+15t+10a\le 11n+12m+3a+1\le 12n+36m+18t+12a$ so that $t\le (n-18m-7a+1)/15$.
\item Suppose $v=v_{r,1}, r\in[1,a]$ with $k=r$. Therefore, we must have $(2j'-1)r \le 11n+12m+3a+1 \le 2j'r$. Note that $s\ge 2$ if $r=1$, and $s\ge 1$ otherwise.
\item Suppose $v=v_{r,3}$ with $k=a+r$. Therefore, we must have $(2j'-1)(a+r)\le 11n+12m+3a+1 \le 2j'(a+r)$. 
\item Suppose $v = u_{i,1}$, $i\in [1,m]$ with $k=2n+2a+3i-2$. Therefore, we must have $(2j'-1)(2n+2a+3i-2)\le 11n+12m+3a+1\le 2j'(2n+3a+3i-2)$. 
\item Suppose $v= u_{i,3}$ with $k=2n+2a+3i-1$. Therefore, we must have $(2j'-1)(2n+2a+3i-1)\le 11n+12m+3a+1 \le 2j'(2n+3a+3i-1)$. 
\item Suppose $v= u_{i,5}$ with $k=2n+2a+3i$. Therefore, we must have $(2j'-1)(2n+2a+3i)\le 11n+12m+3a+1 \le 2j'(2n+3a+3i)$. 
\item Suppose $v= u_{i,2}$ with $k=8n+12m+2a-3i+3$. Therefore, we must have $(2j'-1)(8n+12m+2a-3i+3)\le 11n+12m+3a+1 \le 2j'(8n+12m+2a-3i+3)$. Since $i\le m$, we must have $j'=1$ and the equality always hold. 
\item Suppose $v=u_{i,4}$ with $k=8n+12m+2a-3i+1$. Therefore, we must have $(2j'-1)(8n+12m+2a-3i+1) \le 11n+12m+3a+1 \le 2j'(8n+12m+2a-3i+1)$. Similar to (8), the equality always hold. 
\item Suppose $v=u_{i,6}$ with $k=8n+12m+2a-3i+2$.  Therefore, we must have $(2j'-1)(8n+12m+2a-3i+2)\le 11n+12m+3a+1 \le 2j'(8n+12m+2a-3i+2)$. Similar to (8), the equality always hold. 
\item Suppose $v=y_{t,3}, t\in [1,n]$ with $k=8n+6m+2a-3t+3$. Therefore, we must have $(2j'-1)(8n+6m+2a-3t+3)\le 11n+12m+3a+1 \le 2j'(8n+6m+2a-3t+3)$. Thus, $j'=1$ so that $t\le (5n+a+5)/6$. Since $t\le n$, we have $a+5\le n$. The conditions $n\ge 2a+2$ or else $a\ge 2n$ further implies that we must have $n\ge a+5$ for $a=1,2$, and $n\ge 2a+2$ otherwise. 
\item Suppose $v=y_{t,5}$ with $k=8n+6m+2a-3t+1$. Therefore, we must have $(2j'-1)(8n+6m+2a-3t+1)\le 11n+12m+3a+1 \le 2j'(8n+6m+2a-3t+1)$. Thus, $j'=1$ so that $t\le (5n+a+1)/6$. Since $t\le n$, we must have $a+1\le n$. Similar to (11), we must have $n\ge 2a+2$ for this case.
\item Suppose $v=y_{t,2}$ with $k=2n+6m+2a+3t-2$. Therefore, we must have $(2j'-1)(2n+6m+2a+3t-2)\le 11n+12m+3a+1 \le 2j'(2n+6m+2a+3t-2)$. Thus, $j'\le 3$. \\ If $j'=1$, then $2n+6m+2a+3t-2\le 11n+12m+3a+1\le 4n+12m+4a+6t-4$ so that $(7n-a+5)/6 \le t\le (9n+6m+a+3)/3$. Since $t\le n$, we then have $a\ge n+5$. The conditions $n\ge 2a+2$ or else $a\ge 2n$ further implies that we must have $a\ge n+5$ for $a=1,2,3,4$ and $a\ge 2n$ otherwise. \\ If $s\ge j'=2$, then $6n+18m+6a+9t-6 \le 11n+12m+3a+1\le 8n+24m+8a+12t-8$ so that $(3n-12m-5a+9)/12 \le t \le (5n-6m-3a+7)/9$.  \\ If $s\ge j'=3$, then $10n+30m+10a+15t-10\le 11n+12m+3a+1\le 12n+36m+12a+18t-12$ so that $t\le (n-18m-7a+13)/15$.
\item Suppose $v=y_{t,4}$ with $k=2n+6m+2a+3t-1$. Therefore, we must have $(2j'-1)(2n+6m+2a+3t-1)\le 11n+12m+3a+1 \le 2j'(2n+6m+2a+3t-1)$. Thus, $j'\le 3$. \\ If $j'=1$, then $2n+6m+2a+3t-1\le 11n+12m+3a+1\le 4n+12m+4a+6t-2$ so that $(7n-a+3)/6 \le t\le (9n+6m+a+2)/3$. Since $t\le n$, we then have $a\ge n+3$. The conditions $n\ge 2a+2$ or else $a\ge 2n$ further implies that we must have $a\ge n+3$ for $a=1,2$ and $a\ge 2n$ otherwise. \\ If $s\ge j'=2$, then $6n+18m+6a+9t-3\le 11n+12m+3a+1\le 8n+24m+8a+12t-4$ so that $(3n-12m-5a+5)/12\le t\le (5n-6m-3a+4)/9$. \\ If $s\ge j'=3$, then $10n + 30m + 10a + 15t - 5 \le 11n+12m+3a+1\le 12n+36m+12a+18t-6$ so that $t\le (n-18m-7a+6)/15$.
\item Suppose $v=v_{r,2}, r\in[1,a]$ with $k=12m+11n+3a+1-r$. Therefore, we must have $(2j'-1)(12m+11n+3a+1-r)\le  11n+12m+3a+1 \le 2j'(12m+11n+3a+1-r)$. Since $r\le a$, we must have $j'=1$ and the equality always hold.  
\end{enumerate}
\end{proof}

\begin{example} The figure below gives the local total antimagic 39-coloring of $G_v(k,s)$ for $G=C_6 + P_6 + 3P_3$, $v=v_{3,2}, k=30, s=1$ with induced weights $28,29,33,34,35,36,37,38,69,70,\ldots,98$ and $2576$ being the weight of $v_{3,2}$ in $G_v(k,s)$ as defined above.

\begin{figure}[H]
\begin{center}
\centerline{\epsfig{file=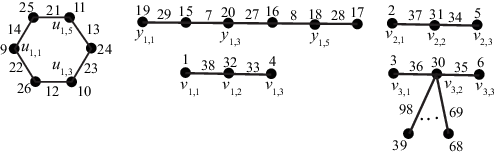, width=12cm}}
\caption{$G_v(k,s)$, $v=v_{3,2}$, $k=30,s=1$ with local total antimatic 39-coloring.}\label{fig:C6P63P3+ks} 
\end{center}
\end{figure}
\end{example}

\section{Conclusions}\label{sec-con}

\nt In this paper, we successfully characterize graphs with $\chi_{lt}(G) = 3$. Moreover, we obtained sufficient condition for a graph $G$ with $k\ge 2$ pendant edges such that $\chi_{lt}(G)\ge k+2$. A family of graph $G$ with $\chi_{lt}(G)\ge k+2$ is given. We then obtained many families of graphs with (i) $\Delta(G)=k+2$ and $\chi_{lt}(G) = k+3$,  or (ii) $k \ge \D(G)$ and $\chi_{lt}(G) = k+1$. The following problems arise naturally.

\begin{problem} For $2n\ge k+3\ge 4$ and $2n(2n+1) - 2n(k+2)(5k+5) + (k+2)(k-1) > 0$, show that $\chi_{lt}(f_n(k)) = 2nk+2$. Otherwise, $\chi_{lt}(f_n(k)) = 2nk+1$. \end{problem}

\begin{problem} Determine $\chi_{lt}(G)$ for $G\cong P_n, C_n$ for $n\ge 3$ completely. \end{problem}

\begin{problem} Determine $\chi_{lt}(G_v(k,s))$ for $G=mC_6+nP_6$, $m\ge 1, n\ge 1$.  \end{problem}

\begin{problem} Determine $\chi_{lt}(mC_6+nP_6+aP_3)$ for $m, a\ge 1, (a+1)/2\le n \le 2a+1$. \end{problem}

\nt  We note that $P_n, n\ge 3$ is the only graph with $\Delta(P_n) = 2 =$ the number of pendant edges. Thus, we also pose the following problem.

\begin{problem} Determine $\chi_{lt}(G)$ for $G$ with $k\ge 3$ pendant edges and $\Delta(G)=k$. \end{problem}

\ms\nt In all the known results on graphs $G$ with $k\ge \D(G)$ pendant edges, we have $k+1\le \chi_{lt}(G)\le k+2$. Similar to the conjecture of $\Delta(G)+1\le \chi_t(G)\le \Delta(G)+2$, we end the paper with the following conjectures.

\begin{conjecture} If $G$ has $k\le \Delta(G)$ pendant edges, then $\Delta(G)+1\le \chi_{lt}(G)\le \Delta(G)+2$. \end{conjecture}

\begin{conjecture} If $G$ has $k\ge \D(G)$ pendant edges, then $k+1\le \chi_{lt}(G)\le k+2$. \end{conjecture}

\ms\ms
\nt{\bf Acknowledgements} 

\ms\nt Most results of this paper were obtained while the author was affiliated to Universiti Teknologi MARA (Segamat Campus). He is grateful to the university for all the supports given.

\end{document}